\documentclass[12pt]{article}

\usepackage{amsmath}
\usepackage{amssymb}
\usepackage{exscale}
\usepackage{amsthm}
\usepackage{epsfig}
\usepackage{verbatim}
\usepackage{setspace}

\DeclareMathOperator{\disc}{disc}

\DeclareMathOperator{\lcm}{lcm}

\theoremstyle{plain}

\newtheorem{theo}{Theorem}

\newtheorem{lemma}[theo]{Lemma}

\newcommand{\mmod}[1]{\,\,(\mbox{mod}\,\,#1)}
\newcommand{\eins}{1\hspace{-1.4mm}1}

\begin{document}

\title{Discrepancy of Sums of two Arithmetic Progressions}

\author{Nils Hebbinghaus\\
Max-Planck-Institut f\"ur Informatik,\\ Saarbr\"ucken, Germany.}
 \date{} 
\maketitle

\begin{abstract}
  Estimating the discrepancy of the hypergraph of all arithmetic
  progressions in the set $[N]=\{1,2,\hdots,N\}$ was one of the famous
  open problems in combinatorial discrepancy theory for a long time.
  An extension of this classical hypergraph is the hypergraph of sums
  of $k$ ($k\geq 1$ fixed) arithmetic progressions. The hyperedges of
  this hypergraph are of the form $A_{1}+A_{2}+\hdots+A_{k}$ in $[N]$,
  where the $A_{i}$ are arithmetic progressions. For this hypergraph
  Hebbinghaus (2004) proved a lower bound of $\Omega(N^{k/(2k+2)})$.
  Note that the probabilistic method gives an upper bound of order
  $O((N\log N)^{1/2})$ for all fixed~$k$. P\v{r}\'{i}v\v{e}tiv\'{y}
  improved the lower bound for all $k\geq 3$ to $\Omega(N^{1/2})$ in
  2005.  Thus, the case $k=2$ (hypergraph of sums of two arithmetic
  progressions) remained the only case with a large gap between the
  known upper and lower bound. We bridge this gap (up to a logarithmic
  factor) by proving a lower bound of order $\Omega(N^{1/2})$ for the
  discrepancy of the hypergraph of sums of two arithmetic
  progressions.
\end{abstract}

\section{Introduction}
A finite hypergraph ${\mathcal H}=(V,{\mathcal E})$ consists of a
finite set $V$ and a set ${\mathcal E}$ of subsets of $V$. The
elements of $V$ are called vertices and those of ${\mathcal E}$
hyperedges of the hypergraph ${\mathcal H}$. If we $2$--partition the
set of vertices $V$, this $2$--partition clearly induces a
$2$--partition in every hyperedge $E\in{\mathcal E}$. The discrepancy
of ${\mathcal H}$ is a non-negative integer that indicates how
ballanced ${\mathcal H}$ can be $2$--partitioned with respect to all
its hyperedges $E\in{\mathcal E}$. We make this more precise and
express the $2$--partition through a coloring $\chi\colon
V\to\{-1,1\}$ of the vertices of ${\mathcal H}$ with the two
``colors'' $-1$ and $1$. Now for every hyperedge $E\in{\mathcal E}$
the imbalance due to the coloring $\chi$ can be calculated as follows.
Let $\chi(E):=\sum_{x\in E}\chi(x)$. Then $|\chi(E)|$ is the absolute
difference between the number of vertices in $E$ colored with
``color'' $-1$ and the number of vertices in $E$ colored with
``color'' $1$. The discrepancy of ${\mathcal H}$ with respect to the
(specific) coloring $\chi$ is defined as
$$\disc({\mathcal H},\chi):=\max_{E\in{\mathcal E}}|\chi(E)|.$$
In other words, $\disc({\mathcal H},\chi)$ is the maximal imbalance of
any hyperedge $E\in{\mathcal E}$ under the coloring $\chi$. Now the
discrepancy of ${\mathcal H}$ is defined as
$$\disc({\mathcal H}):=\min_{\chi\colon V\to\{-1,1\}}\disc({\mathcal H},\chi),$$
where the minimum is taken over all $2^{|V|}$ possible colorings
$\chi\colon V\to\{-1,1\}$ of the set of vertices $V$. Thus, $\disc({\mathcal
  H})$ is the least possible imbalance of any hyperedge $E\in{\mathcal
  E}$ that can not be avoided under any coloring $\chi\colon V\to\{-1,1\}$.

One of the famoust long-standing open problems in (combinatorial)
discrepancy theory was to determine the right order for the
discrepancy of the hypergraph of arithmetic progressions in the first
$N$ natural numbers ($N\in\mathbb{N}$). Before we give a brief
overview over the history of this problem, we introduce the hypergraph
${\mathcal H}_{AP}$ of arithmetic progressions. For convenience, let
us define for every interval $I\subseteq\mathbb{R}$ the set
$$I_{\mathbb{Z}}:=\{z\in\mathbb{Z}\mid z\in I\}$$
of all integers in the intervall $I$. In particular, we introduce the
abreviation
$$[x]:=[1,x]_{\mathbb{Z}}$$
for the set of all natural numbers $n$ with $1\leq n\leq x$
($x\in\mathbb{R}$). Let $N\in\mathbb{N}$. An arithmetic progression in
$[N]$ is a subset of $[N]$ of the form
$$A_{a,\delta,L}:=\{a+j\delta\mid j\in[0,L-1]_{\mathbb{Z}}\}.$$
Now we can define the hypergraph ${\mathcal H}_{AP}=([N],{\mathcal
  E}_{AP})$. The set of vertices of ${\mathcal H}_{AP}$ is the set
$[N]$ and the set of hyperedges is
$${\mathcal E}_{AP}:=\{A_{a,\delta,L}\mid a,\delta\in[N],L\in[\tfrac{N-a}{\delta}+1]\},$$
where $L\in[\tfrac{N-a}{\delta}+1]$ just ensures that
$A_{a,\delta,L}\subseteq [N]$.

In 1964, Roth~\cite{roth} proved a lower bound for the
discrepancy of the hypergraph ${\mathcal H}_{AP}$ of order
$\Omega(N^{\frac{1}{4}})$. Using a random coloring of the vertices of
${\mathcal H}_{AP}$, one can easily show an upper bound of order
$O((N\log N)^{\frac{1}{2}})$ for the discrepancy of the hypergraph
${\mathcal H}_{AP}$. The first non-trivial upper bound is due to
S\'{a}rk\"{o}zy. 1973 he proved $\disc({\mathcal
  H}_{AP})=O(N^{1/3}\log^{1/3} N)$ A sketch of his beautiful proof can
be found in the book {\em Probabilistic Methods in Combinatorics} by
Erd\H{o}s and Spencer~\cite{es74}. Beck~\cite{beck} showed
in 1981 (inventing the famous partial coloring method) that
Roth's lower bound is almost sharp. His upper bound of order
$O(N^{1/4}\log^{5/4} N)$ was finally improved by Matou\v{s}ek and
Spencer~\cite{ms96} in 1996.  They showed by a refinement of
the partial coloring method --- the entropy method --- that the
discrepancy of the hypergraph ${\mathcal H}_{AP}$ is exactly of order
$\Theta(N^{1/4/})$.

Therefore, after 32 years, this open problem was solved. In the next
years several extensions of this discrepancy problem were studied.
Doerr, Srivastav and Wehr~\cite{dsw} determined the discrepancy of
$d$--dimensional arithmetic progressions. For the hypergraph
${\mathcal H}_{AP,d}=([N]^{d},{\mathcal E}_{AP,d})$, where ${\mathcal
  E}_{AP,d}:=\{\prod_{i=1}^{d} E_{i}\mid E_{i}\in{\mathcal E}_{AP}\}$,
they proved $\disc({\mathcal H}_{AP,d})=\Theta(N^{d/4})$. Another
related hypergraph --- the hypergraph of all $1$--dimensional
arithmetic progressions in the $d$--dimensional grid $[N]^{d}$ was
studied by Valko~\cite{v02}. He proved for the discrepancy of this
hypergraph a lower bound of order $\Omega(N^{d/(2d+2)})$ and an upper
bound of order $O(N^{d/(2d+2)}\log^{5/2}N)$.

The hypergraph that we consider in this paper was introduced by
Hebbinghaus~\cite{h04} in a generalized version. Let $k\in\mathbb{N}$
and $N\in\mathbb{N}$. The hypergraph ${\mathcal
  H}_{kAP}=([N],{\mathcal E}_{kAP})$ of sums of $k$ arithmetic
progressions is defined as follows. The vertices of ${\mathcal
  H}_{kAP}$ are the first $N$ natural numbers. And the set of
hyperedges ${\mathcal E}_{kAP}$ is defined as
$${\mathcal E}_{kAP}:=\left\{\left(\left.\sum_{i=1}^{k}A_{a_{i},\delta_{i},L_{i}}\right)\cup[N]\right| a_{i}\in\mathbb{Z},\delta_{i},L_{i}\in[N](i\in[k])\right\},$$
where the sum of $k$ sets $M_{i}$ ($i\in[k]$) is
$$\sum_{i=1}^{k}M_{i}=\left\{\left.\sum_{i=1}^{k}\right| m_{i}\in M_{i}(i\in[k])\right\}.$$
For the hypergraph ${\mathcal H}_{kAP}$ of sums of $k$ arithmetic
progressions in $[N]$ Hebbinghaus~\cite{h04} proved a lower bound of
order $\omega(N^{k/(2k+2)})$ in 2004. But there remained a large gap
between this bound and the upper bound of order
$O(N^{1/2}\log^{1/2}N)$ from the random coloring method. In 2006
P\v{r}\'{i}v\v{e}tiv\'{y}~\cite{p06} nearly closed this gap for $k\geq
3$ by proving a lower bound of order $\Omega(N^{1/2})$ for the
discrepancy of the hypergraph ${\mathcal H}_{3AP}$ of sums of three
arithmetic progressions. This lower bound clearly extends to all
hypergraphs ${\mathcal H}_{kAP}$ for all $k\geq 3$. Thus, the case
$k=2$ was the last with a large gap between the lower and the upper
bound for the discrepancy. In this paper we improve the lower bound
for the discrepancy of the hypergraph ${\mathcal H}_{2AP}$ of sums of
two arithmetic progressions from the order $\Omega(N^{1/3})$ to the
order $\Omega(N^{1/2})$. This result shows that the upper bound of
order $O(N^{1/2}\log^{1/2}N)$ for the discrepancy of ${\mathcal
  H}_{2AP}$ determined by the random coloring method is almost sharp.
We will prove the following theorem.
\begin{theo}\label{MainTheo}
  Let $N\in\mathbb{N}$. For the hypergraph ${\mathcal H}_{2AP}$ of
  sums of two arithmetic progressions we obtain the following bounds.
  \begin{itemize}
  \item[(i)] $\disc({\mathcal H}_{2AP})=\Omega(N^{1/2})$.
  \item[(ii)] $\disc({\mathcal H}_{2AP})=O(N^{1/2}\log^{1/2}N)$.
  \end{itemize}
\end{theo}
Since $|{\mathcal E}_{2AP}|=O(N^{6})$ the second assertion is a
direct consequence of the general upper bound for a hypergraph
${\mathcal H}$ with $n$ vertices and $m$ hyperedges $\disc({\mathcal
  H})=O(\sqrt{n\log m})$ derived by the random coloring method.

\section{A Special Set of Hyperedges}\label{sec_E0}
In this section we define a special subset ${\mathcal E}_{0}$ of the
set ${\mathcal E}_{2AP}$ of all sums of two arithmetic progressions in
$[N]$. The elements of this set ${\mathcal E}_{0}$ and all their
translates build the set of hyperedges in which we will find for every
coloring $\chi\colon  V\to\{-1,1\}$ a hyperedge with discrepancy of order
$\Omega(N^{\frac{1}{2}})$.

All elements of ${\mathcal E}_{0}$ are sums of two arithmetic
progressions with starting point $0$. Thus, we can characterize them
by the difference and length of the two arithmetic progressions. We
define for all $\delta_{1}, \delta_{2}, L_{1}, L_{2}\in\mathbb{N}$:
$$E_{\delta_{1}, L_{1}, \delta_{2}, L_{2}}:=\{j_{1}\delta_{1}+j_{2}\delta_{2}\mid j_{1}\in[0,L_{1}-1], j_{2}\in[0,L_{2}-1]\}.$$
Before specifying the set ${\mathcal E}_{0}$ we should mention that
due to a case distinction in the proof of the Main Lemma the set
${\mathcal E}_{0}$ is the union of three subsets ${\mathcal E}_{1}$,
${\mathcal E}_{2}$ and ${\mathcal E}_{3}$, each of them corresponding
to one of the cases. The first two sets ${\mathcal E}_{1}$ and
${\mathcal E}_{2}$ are easy to define. We set
$${\mathcal E}_{1}:=\{E_{\delta_{1}, L_{1}, \delta_{2}, L_{2}}\mid \delta_{1}\in[24], L_{1}=\lceil\tfrac{N}{6\delta_{1}}\rceil,\delta_{2}=1, L_{2}=1\},$$
and
$${\mathcal E}_{2}:=\{E_{\delta_{1}, L_{1}, \delta_{2}, L_{2}}\mid \delta_{1}\in[25,N^{\frac{1}{2}}]_{\mathbb{Z}}, L_{1}=\lceil\tfrac{N}{12\delta_{1}}\rceil,\delta_{2}\in[\delta_{1}-1], L_{2}=\lceil\tfrac{\delta_{1}-1}{12}\rceil\}.$$
The definition of the last set ${\mathcal E}_{3}$ is not
straightforward. For every difference $\delta_{1}$ of the first
arithmetic progression we have to determine a set of differences
$\delta_{2}$ for the second arithmetic progression. Let
$\delta_{1}\in[N^{\frac{1}{2}}]$ and let
$$B(\delta_{1}):=\{b\in[\delta_{1}]\mid (b,\delta_{1})=1\}$$
be the set of all elements of $[\delta_{1}]$ that are relatively prime
to $\delta_{1}$. Here $(b,\delta_{1})$ denotes the greatest common
divisor of $b$ and $\delta_{1}$. Let $b\in B(\delta_{1})$. Set
$\bar{k}:=\lfloor\log(N^{\frac{1}{2}}\delta_{1}^{-1})\rfloor$. We
define for all $0\leq k\leq \bar{k}$ sets $M(b,k)$ of distances for
the second arithmetic progression. The set $M(b,k)$ should cover the
range of possible differences for the second arithmetic progression
for the interval $(2^{k}N^{\frac{1}{2}},2^{k+1}N^{\frac{1}{2}}]$. We
define
$$M(b,k):=(b+2^{2k}\delta_{1}\mathbb{Z})\cap(2^{k}N^{\frac{1}{2}},2^{k+1}N^{\frac{1}{2}}+2^{2k}\delta_{1}).$$
For all $0\leq k\leq\bar{k}$, we set
$M_{\delta_{1}}(k):=\bigcup\limits_{b\in B(\delta_{1})}M(b,k)$. Now we
are able to define the third set ${\mathcal E}_{3}$. Let
$${\mathcal E}_{3}:=\bigcup_{\delta_{1}\in[N^{\frac{1}{2}}]}\bigcup_{k=0}^{\bar{k}}\left\{E_{\delta_{1},L_{1},\delta_{2},L_{2}}\mid L_{1}=\left\lceil\tfrac{2^{k}N^{\frac{1}{2}}}{12}\right\rceil,\delta_{2}\in M_{\delta_{1}}(k),L_{2}=\left\lceil\tfrac{2^{-k}N^{\frac{1}{2}}}{12}\right\rceil\right\}.$$
In the next lemma we prove that the cardinality of the set ${\mathcal
  E}_{0}$ is of order $O(N)$. This is an essential property of the set
${\mathcal E}_{0}$ for the proof of the lower bound of the discrepancy
of the hypergraph of sums of two arithmetic progressions.
\begin{lemma}\label{propertiesE3}
  We have $|{\mathcal E}_{3}|\leq 6N$ and thus
  $|{\mathcal E}_{0}|\leq 7N$.
\end{lemma}
\begin{proof}
  We have to estimate $|{\mathcal
    E}_{3}|=\sum_{\delta_{1}\in\left[N^{\frac{1}{2}}\right]}\sum_{k=0}^{\bar{k}}|M_{\delta_{1}}(k)|$.
  For this pupose we look for $|M(b,k)|$ for all $b\in B(\delta_{1})$
  and all $0\leq k\leq\bar{k}$. We first show that the difference
  $2^{2k}\delta_{1}$ of two consecutive elements of $M(b,k)$ is at
  most $2^{k}N^{\frac{1}{2}}$.
  $$2^{2k}\delta_{1}\leq 2^{k}2^{\log\left(N^{\frac{1}{2}}\delta_{1}^{-1}\right)}\delta_{1}=2^{k}N^{\frac{1}{2}}.$$
  Hence,
  $$|M(b,k)|\leq\frac{3\cdot 2^{k}N^{\frac{1}{2}}}{2^{2k}\delta_{1}}=3\cdot 2^{-k}N^{\frac{1}{2}}\delta_{1}^{-1}.$$
  Since $M_{\delta_{1}}(k)=\bigcup\limits_{b\in B(\delta_{1})}M(b,k)$, this yields $|M_{\delta_{1}}(k)|\leq\delta_{1}|M(b,k)|\leq 3\cdot 2^{-k}N^{\frac{1}{2}}$. Thus, we get
  \begin{eqnarray*}
    |{\mathcal E}_{3}|&=&\sum_{\delta_{1}\in\left[N^{\frac{1}{2}}\right]}\sum_{k=0}^{\bar{k}}|M_{\delta_{1}}(k)|\\
    &\leq&\sum_{\delta_{1}\in\left[N^{\frac{1}{2}}\right]}\sum_{k=0}^{\bar{k}}3\cdot 2^{-k}N^{\frac{1}{2}}\\
    &<&3N\sum_{k=0}^{\infty} 2^{-k}\\
    &\leq&6N
  \end{eqnarray*}
  It is easy to see that $|{\mathcal E}_{1}\cup{\mathcal E}_{2}|<N$.
  This proves the lemma.
\end{proof}

\section{Discrete Fourier Analysis}
The purpose of this section is discrete Fourier analysis on the
additive group $(\mathbb{Z},+)$ and its connection to the discrepancy
of the hypergraph ${\mathcal H}_{2AP}$. First of all, let us extend
the coloring $\chi$ to the set of all integers as follows. We keep the
``old'' color values for the set $[N]$ and set $\chi(z):=0$ for all
$z\in\mathbb{Z}\setminus[N]$. Thus, the (extended) coloring
$\chi\colon \mathbb{Z}\to\{-1,0,1\}$ satisfies the condition:
$\chi(z)=0$, if and only if $z\in\mathbb{Z}\setminus[N]$. For every
set $E\subseteq\mathbb{Z}$ we define its color value
$\chi(E):=\sum_{x\in E}\chi(x)$. One can easily verify that we can
express the coloring value of the set $E_{a}:=a+E=\{a+x\mid x\in E\}$
as convolution of~$\chi$ and the indicator function $\eins_{-E}$ of
the set $-E=\{-x\mid x\in E\}$ evaluated at $a$.  For all
$a\in\mathbb{Z}$, we have
\begin{eqnarray}\label{convolution}
        \chi(E_{a})=(\chi\ast\eins_{-E})(a).
\end{eqnarray}
Thus, for all $E\subseteq\mathbb{Z}$
$$\sum_{a\in\mathbb{Z}}|\chi(E_{a})|^{2}=\|\chi\ast\eins_{-E}\|_{2}^{2}.$$
For the proof of the lower bound in Theorem~\ref{MainTheo} we use a $2$--norm
approach. More precisely, we will estimate the sum of squarred
discrepancies
\begin{eqnarray}\label{sumsquarreddisc}
        \sum_{E\in{\mathcal
E}_{0}}\sum_{a\in\mathbb{Z}}|\chi(E_{a})|^{2}=\sum_{E\in{\mathcal
E}_{0}}\|\chi\ast\eins_{-E}\|_{2}^{2}.
\end{eqnarray}
Using two well-known facts from Fourier analysis, the Plancherel
Theorem and the multiplicity of the Fourier transform, we will lower
bound this sum of squarred discrepancies. Afterwards an averaging
argument will yield the existence of a hyperedge $E$ with a
discrepancy of order $\Omega(N^{\frac{1}{2}})$. But first of all we
introduce the Fourier transform of a function $f\colon
\mathbb{Z}\to\mathbb{C}$. The Fourier transform of $f$ is defined as
$$\widehat{f}\colon [0,1)\to\mathbb{C},\quad
\alpha\mapsto\sum_{z\in\mathbb{Z}}f(z)e^{2\pi iz\alpha}.$$
In the following lemma we list the two facts from Fourier analysis on the
additive group $(\mathbb{Z},+)$ that we will need for our calculations.
\begin{lemma}\label{LemmaFourier}
  Let $f,g\colon \mathbb{Z}\to\mathbb{C}$ two square integrable
  functions. For the Fourier transform of $f$ and $g$ we get
        \begin{itemize}
        \item[(i)] $\|\widehat{f}\|_{2}^{2}=\|f\|_{2}^{2}$ (Plancherel
          Theorem),
                \item[(ii)] $\widehat{f\ast g}=\widehat{f}\widehat{g}$.
        \end{itemize}
\end{lemma}

\section{Proof of the Lower Bound}
Before we prove the lower bound for the discrepancy of the hypergraph
of sums of two arithmetic progressions, we state the following lemma.
\begin{lemma}[Main Lemma]\label{mainlemma}
  For every $\alpha\in[0,1)$, there exists an $E\in{\mathcal E}_{0}$
  such that
  $$|\widehat{\eins}_{-E}(\alpha)|\geq\frac{1}{300}N.$$
\end{lemma}
Applying this lemma, we are able to give the lower bound proof.
\begin{proof}[Proof of Theorem~\ref{MainTheo}] Using the
equation~(\ref{sumsquarreddisc}) and Lemma~\ref{LemmaFourier}, we get
\begin{eqnarray*}
  \sum_{E\in{\mathcal
      E}_{0}}\sum_{a\in\mathbb{Z}}|\chi(E_{a})|^{2}&=&\sum_{E\in{\mathcal
      E}_{0}}\|\chi\ast\eins_{-E}\|_{2}^{2}\\
  &=&\sum_{E\in{\mathcal E}_{0}}\|\widehat{\chi\ast\eins}_{-E}\|_{2}^{2}\\
  &=&\sum_{E\in{\mathcal
      E}_{0}}\|\widehat{\chi}\widehat{\eins}_{-E}\|_{2}^{2}\\
  &=&\sum_{E\in{\mathcal
      E}_{0}}\int_{0}^{1}|\widehat{\chi}(\alpha)|^{2}|\widehat{\eins}_{-E}(\alpha)|^{2}d\alpha\\
  &=&\int_{0}^{1}|\widehat{\chi}(\alpha)|^{2}\left(\sum_{E\in{\mathcal
        E}_{0}}|\widehat{\eins}_{-E}(\alpha)|^{2}\right)d\alpha.
\end{eqnarray*}
The Main Lemma yields for every $\alpha\in[0,1)$ the existence of an
$E\in{\mathcal E}_{0}$ such that
$|\widehat{\eins}_{-E}(\alpha)|\geq\tfrac{1}{300}N$. Thus, we get for
every $\alpha\in[0,1)$
$$\sum_{E\in{\mathcal
E}_{0}}|\widehat{\eins}_{-E}(\alpha)|^{2}\geq\tfrac{1}{90000}N^{2}.$$
Hence, we can continue the estimation of the sum of squarred discrepancies
as follows.
\begin{eqnarray*}
        \sum_{E\in{\mathcal E}_{0}}\sum_{a\in\mathbb{Z}}|\chi(E_{a})|^{2}
        &=&\int_{0}^{1}|\widehat{\chi}(\alpha)|^{2}\left(\sum_{E\in{\mathcal
E}_{0}}|\widehat{\eins}_{-E}(\alpha)|^{2}\right)d\alpha\\
        &\geq&\tfrac{1}{90000}N^{2}\|\widehat{\chi}\|_{2}^{2}\\
        &=&\tfrac{1}{90000}N^{2}\|\chi\|_{2}^{2}\\
        &=&\tfrac{1}{90000}N^{3}.
\end{eqnarray*}
Since every $E\in{\mathcal E}_{0}$ satisfies
$E\subseteq[0,N-1]_{\mathbb{Z}}$, we get for every
$a\in\mathbb{Z}\setminus[-N+1,N]_{\mathbb{Z}}$ that $E\cap [N]=\emptyset$
and thus $\chi(E_{a})=0$. Therefore, $\sum_{E\in{\mathcal
E}_{0}}\sum_{a\in\mathbb{Z}}|\chi(E_{a})|^{2}$ is the sum of at most
$2N|{\mathcal E}_{0}|\leq 14 N^{2}$ non-trivial elements
(Lemma~\ref{propertiesE3}). Hence, there exists an $E\in{\mathcal E}_{0}$ and an
$a\in[-N+1,N]_{\mathbb{Z}}$ such that
$$|\chi(E_{a})|^{2}\geq\tfrac{1}{1260000}N.$$
Thus, we have proven
$$\disc({\mathcal
H}_{2AP})\geq|\chi(E_{a})|>\tfrac{1}{1200}N^{\frac{1}{2}}.$$
\end{proof}

Before we can prove the Main Lemma, we have to state
and prove the following four lemmas.
\begin{lemma}\label{pigeonhole}
  For every $\alpha\in[0,1)$ and every $k\in\mathbb{N}$, there exists a
  $\delta\in[k]$ and an $a\in\mathbb{Z}$ such that
  $$|\delta\alpha-a|<\frac{1}{k}.$$
\end{lemma}
\begin{proof}
  For all $j\in[k]$, we define
  $$M_{j}:=\left\{\delta\in[k]\colon \delta\alpha-\lfloor\delta\alpha\rfloor\in\left[\frac{j-1}{k},\frac{j}{k}\right)\right\}.$$
  For every $\delta\in M_{1}$, holds
  $|\delta\alpha-\lfloor\delta\alpha\rfloor|<\tfrac{1}{k}$. Thus, we
  can assume $M_{1}=\emptyset$. By the pigeon hole principle, there
  exists a $j\in[k]\setminus\{1\}$ with $|M_{j}|\geq 2$. Let
  $\delta_{1},\delta_{2}\in M_{j}$ with $\delta_{1}<\delta_{2}$. Set
  $\delta:=\delta_{2}-\delta_{1}$. Using $\delta_{1},\delta_{2}\in
  M_{j}$, we get
  $$|\delta-(\lfloor\delta_{2}\alpha\rfloor-\lfloor\delta_{1}\alpha\rfloor)|=|(\delta_{2}-\lfloor\delta_{2}\alpha\rfloor)-(\delta_{1}-\lfloor\delta_{1}\alpha\rfloor)|<\frac{1}{k}.$$
\end{proof}
\begin{lemma}\label{modulo}
  Let $a,\delta\in\mathbb{N}$ with $(a,\delta)=1$. There exists a
  $k\in[\delta-1]$ such that
  $$ka\equiv 1 \mmod{\delta}.$$
  Moreover, $(\delta-k)a\equiv -1\mmod{\delta}$. It holds
  $(k,\delta)=(\delta-k,\delta)=1$.
\end{lemma}
\begin{proof}
  Since $(a,\delta)=1$, there exist $k,\ell\in\mathbb{Z}$ with 
  $$ka+\ell\delta=1.$$
  Thus, $ka\equiv 1 \mmod{\delta}$. Obviously, $k$ can be choosen from
  the set $[\delta-1]$. The second assertion follows from
  $$ka+(\delta-k)a=\delta a\equiv 0 \mmod{\delta}.$$
  Finally, the equation $ka+\ell\delta=1$ proves also $(k,\delta)=1$. But
  this implies $(\delta-k,\delta)=1$.
\end{proof}
\begin{lemma}\label{fourier}
  Let $\alpha\in[0,1)$,
  $\delta_{1},\delta_{2},L_{1},L_{2}\in\mathbb{N}$ with $L_{1}\neq
  1\neq L_{2}$ be chosen such that for suitable
  $a_{1},a_{2}\in\mathbb{Z}$ we have
  $$|\delta_{j}\alpha-a_{j}|\leq\frac{1}{12(L_{j}-1)},\quad (j=1,2).$$
  Set
  $E:=\{j_{1}\delta_{1}+j_{2}\delta_{2}\colon j_{1}\in[0,L_{1}-1],j_{2}\in[0,L_{2}-1]\}$.
  For the Fourier transform of the indicator function $\eins_{-E}$ of
  the set $-E$ we get
  $$|\widehat{\eins}_{-E}(\alpha)|\geq\frac{|E|}{2}.$$
\end{lemma}
\begin{proof}
  The Fourier transform of a function $f\colon \mathbb{Z}\to\mathbb{C}$ is
  given as $\widehat{f}\colon [0,1)\to\mathbb{C}$,
  $\alpha\mapsto\sum\limits_{z\in\mathbb{Z}}f(z)e^{-2\pi i z\alpha}$. Thus,
  $$\widehat{\eins}_{-E}(\alpha)=\sum\limits_{z\in E}e^{2\pi
    iz\alpha}.$$
  Let $z\in E$. There exists a $j_{1}\in[0,L_{1}-1]$ and a
  $j_{2}\in[0,L_{2}-1]$ with
  $z=j_{1}\delta_{1}+j_{2}\delta_{2}$. Hence,
  \begin{eqnarray*}
    e^{2\pi iz\alpha}&=&e^{2\pi i(j_{1}\delta_{1}+j_{2}\delta_{2})\alpha}\\
    &=&e^{2\pi i[j_{1}(\delta_{1}\alpha-a_{1})+j_{2}(\delta_{2}\alpha-a_{2})]}e^{2\pi i(j_{1}a_{1}+j_{2}a_{2})}\\
    &=&e^{2\pi i[j_{1}(\delta_{1}\alpha-a_{1})+j_{2}(\delta_{2}\alpha-a_{2})]}.
  \end{eqnarray*}
  Using
  $|j_{1}(\delta_{1}\alpha-a_{1})+j_{2}(\delta_{2}\alpha-a_{2})|\leq
  \tfrac{L_{1}-1}{12(L_{1}-1)}+\tfrac{L_{2}-1}{12(L_{2}-1)}=\tfrac{1}{12}+\tfrac{1}{12}=\tfrac{1}{6}$,
  we get $\Re(e^{2\pi iz\alpha})\geq\tfrac{1}{2}$. This proves
  $$|\widehat{\eins}_{-E}(\alpha)|=\sum\limits_{z\in E}e^{2\pi iz\alpha}\geq\Re(\sum\limits_{z\in E}e^{2\pi iz\alpha})\geq\frac{|E|}{2}.$$
\end{proof}
\begin{lemma}\label{cardinality}
  Let $\delta_{1},\delta_{2},L_{1},L_{2}\in\mathbb{N}$. If
  $L_{1}\leq\tfrac{\delta_{2}}{(\delta_{1},\delta_{2})}$ then
  $$|\{j_{1}\delta_{1}+j_{2}\delta_{2}\colon j_{1}\in[0,L_{1}-1],j_{2}\in[0,L_{2}-1]\}|=L_{1}L_{2}.$$
\end{lemma}
\begin{proof}
  Assume there are
  $(j_{1},j_{2}),(j'_{1},j'_{2})\in[0,L_{1}-1]\times[0,L_{2}-1]$ such
  that $(j_{1},j_{2})\neq(j'_{1},j'_{2})$ and
  $$j_{1}\delta_{1}+j_{2}\delta_{2}=j'_{1}\delta_{1}+j'_{2}\delta_{2}.$$
  Clearly, $j_{1}\neq j'_{1}$ and $j_{2}\neq j'_{2}$. Since
  $(j_{1}-j'_{1})\delta_{1}=(j'_{2}-j_{2})\delta_{2}$ is divisible by
  $\delta_{1}$ and $\delta_{2}$ and thus also by their least common
  multiple
  $\lcm(\delta_{1},\delta_{2})=\tfrac{\delta_{1}\delta_{2}}{(\delta_{1},\delta_{2})}$,
  we get
  $$L_{1}>|j_{1}-j'_{1}|\geq\frac{\delta_{2}}{(\delta_{1},\delta_{2})}.$$
  This contradiction shows that the function
  $$f\colon [0,L_{1}-1]\times[0,L_{1}-1]\to\mathbb{Z},\quad (j_{1},j_{2})\mapsto
  j_{1}\delta_{1}+j_{2}\delta_{2}$$ is injective which proves the
  assumption.
\end{proof}
By combining Lemma~\ref{pigeonhole}, Lemma~\ref{modulo},
Lemma~\ref{fourier}, and Lemma~\ref{cardinality}, we are able to prove
the Main Lemma. Recall that we proved the lower bound for the
discrepancy of the hypergraph of all sums of two arithmetic
progressions just by applying the Main Lemma.
\begin{proof}[Proof of the Main Lemma.]
  Using Lemma~\ref{pigeonhole}, we can find a
  $\delta_{1}\in[N^{\frac{1}{2}}]$ such that for an appropriate
  $a_{1}\in\mathbb{Z}$ it holds
  $|\delta_{1}\alpha-a_{1}|<N^{-\frac{1}{2}}$. Dividing by
  $\delta_{1}$, we get
  \begin{eqnarray}\label{delta1}
    |\alpha-\frac{a_{1}}{\delta_{1}}|<N^{-\frac{1}{2}}\delta_{1}^{-1}.
  \end{eqnarray}
  We can choose $\delta_{1}$ and $a_{1}$ in such a way that
  $\tfrac{a_{1}}{\delta_{1}}$ is an irreducible fraction. We
  distinguish three cases.

  {\bf Case $1$}: $|\alpha-\tfrac{a_{1}}{\delta_{1}}|<N^{-1}$ and
  $\delta_{1}\leq 24$.

  Set $L_{1}:=\lceil\tfrac{N}{6\delta_{1}}\rceil$, $\delta_{2}:=1$,
  and $L_{2}:=1$. The set $E:=E_{\delta_{1},L_{1},\delta_{2},L_{2}}$
  is an element of the special set of hyperedges ${\mathcal E}_{0}$.
  More precisely, $E\in{\mathcal E}_{1}$. Arguments similar to those
  used in the proof of Lemma~\ref{fourier} show
  \begin{eqnarray*}
    |\widehat{\eins}_{-E}(\alpha)|&\geq&\Re\left(\sum\limits_{z\in
        E}e^{2\pi i z\alpha}\right)\\
    &=&\sum\limits_{j_{1}=0}^{L_{1}-1}\Re\left(e^{2\pi
        ij_{1}\delta_{1}\alpha}\right)\\
    &\geq&L_{1}\Re\left(e^{\frac{2\pi i}{6}}\right)\\
    &\geq&\frac{N}{288}.
  \end{eqnarray*}
  
  {\bf Case $2$}: $|\alpha-\tfrac{a_{1}}{\delta_{1}}|<N^{-1}$ and
  $\delta_{1}>24$.

  Set $L_{1}:=\lceil\tfrac{N}{12\delta_{1}}\rceil$. Using again
  Lemma~\ref{pigeonhole}, there is a $\delta_{2}\in[\delta_{1}-1]$ such
  that for a suitable $a_{2}\in\mathbb{Z}$ it holds
  $|\delta_{2}\alpha-a_{2}|\leq\tfrac{1}{\delta_{1}-1}$. Hence,
  \begin{eqnarray}\label{delta2}
    |\alpha-\frac{a_{2}}{\delta_{2}}|\leq\frac{1}{(\delta_{1}-1)\delta_{2}}.
  \end{eqnarray}
  Set $L_{2}:=\lceil\tfrac{\delta_{1}-1}{12}\rceil$. Since
  $\tfrac{a_{1}}{\delta_{1}}$ is an irreducible fraction and
  $\delta_{2}<\delta_{1}$ we get
  $\tfrac{a_{1}}{\delta_{1}}\neq\tfrac{a_{2}}{\delta_{2}}$. Thus,
  \begin{eqnarray}\label{a1a2gr}
    |\frac{a_{1}}{\delta_{1}}-\frac{a_{2}}{\delta_{2}}|\geq\frac{1}{\lcm(\delta_{1},\delta_{2})}.
  \end{eqnarray}
  On the other hand, using~(\ref{delta1}) and~(\ref{delta2}) we get
  \begin{eqnarray}\label{a1a2kl}
    |\frac{a_{1}}{\delta_{1}}-\frac{a_{2}}{\delta_{2}}|\leq|\frac{a_{1}}{\delta_{1}}-\alpha|+|\alpha-\frac{a_{2}}{\delta_{2}}|&\leq& \frac{1}{N\delta_{1}}+\frac{1}{(\delta_{1}-1)\delta_{2}}\nonumber\\
    &\leq&\left(\frac{1}{\delta_{1}}+\frac{\delta_{1}}{\delta_{1}-1}\right)\frac{1}{\delta_{1}\delta_{2}}\\
    &<&\frac{13}{12}\frac{1}{\delta_{1}\delta_{2}}\nonumber
  \end{eqnarray}
  Combining~(\ref{a1a2gr}) and~(\ref{a1a2kl}) gives
  $$\frac{1}{\lcm(\delta_{1},\delta_{2})}<\frac{13}{12}\frac{1}{\delta_{1}\delta_{2}}.$$
  But this implies
  $(\delta_{1},\delta_{2})=\tfrac{\delta_{1}\delta_{2}}{\lcm(\delta_{1},\delta_{2})}<\tfrac{13}{12}$ and thus $(\delta_{1},\delta_{2})=1$. Define
  $$E:=E_{\delta_{1},L_{1},\delta_{2},L_{2}}=\{j_{1}\delta_{1}+j_{2}\delta_{2}\colon j_{1}\in[0,L_{1}-1],j_{2}\in[0,L_{2}-1]\}.$$
  We have $E\in{\mathcal E}_{2}\subseteq{\mathcal E}_{0}$. Since
  $L_{2}=\lceil\tfrac{\delta_{1}-1}{12}\rceil<\tfrac{\delta_{1}}{6}<\tfrac{\delta_{1}}{(\delta_{1},\delta_{2})}$,
  we can apply Lemma~\ref{cardinality} and get
  $|E|=L_{1}L_{2}\geq\tfrac{N}{12\delta_{1}}\tfrac{\delta_{1}-1}{12}\geq\tfrac{1}{150}N$.
  Furthermore, $\delta_{1},\delta_{2},L_{1},L_{2}$ satisfy the
  conditions of Lemma~\ref{fourier}. Thus,
  $$|\widehat{\eins}_{-E}(\alpha)|\geq\frac{|E|}{2}\geq\frac{1}{300}N.$$

  {\bf Case $3$}: $|\alpha-\tfrac{a_{1}}{\delta_{1}}|\geq N^{-1}$.

  Choose $k$ such that
  \begin{eqnarray}\label{alphak}
    \left|\alpha-\frac{a_{1}}{\delta_{1}}\right|\in[2^{-k-1}N^{-\frac{1}{2}}\delta_{1}^{-1},2^{-k}N^{-\frac{1}{2}}\delta_{1}^{-1}).
  \end{eqnarray}
  Since $|\alpha-\tfrac{a_{1}}{\delta_{1}}|$ is lower bounded by
  $N^{-1}$ and from above by $N^{-\frac{1}{2}}\delta_{1}^{-1}$ (by
  Inequality~(\ref{delta1})), it holds $0\leq k\leq
  \log(N^{\frac{1}{2}}\delta_{1}^{-1})$. Set $L_{1}:=\lceil
  \tfrac{2^{k}N^{\frac{1}{2}}}{12}\rceil$. Using
  $(a_{1},\delta_{1})=1$, we can apply Lemma~\ref{modulo}, which yields
  the existence of a $\gamma\in[\delta_{1}-1]$ such that
  \begin{itemize}
  \item [(i)] $\gamma a_{1}\equiv 1\mmod{\delta_{1}}$,
  \item [(ii)] $(\delta_{1}-\gamma)a_{1}\equiv -1\mmod{\delta_{1}}$.
  \end{itemize}

  Let
  $s:=(\alpha-\tfrac{a_{1}}{\delta_{1}})|\alpha-\tfrac{a_{1}}{\delta_{1}}|^{-1}$,
  i.e., $s$ is the algebraic sign of
  $(\alpha-\tfrac{a_{1}}{\delta_{1}})$. If $s=1$ we set
  $b:=\delta_{1}-\gamma$, otherwise we set $b:=\gamma$. In both cases
  there exists a $\mu\in\mathbb{Z}$ such that
  $$b\frac{a_{1}}{\delta_{1}}=\mu-\frac{s}{\delta_{1}}.$$
  Define
  $d:=|\alpha-\tfrac{a_{1}}{\delta_{1}}|^{-1}2^{-2k}\delta_{1}^{-2}-b
  2^{-2k}\delta_{1}^{-1}$ and $\delta_{2}:=b+\lceil d\rceil
  2^{2k}\delta_{1}$. Then
  \begin{eqnarray*}
    \delta_{2}\alpha&=&(b+\lceil d\rceil
    2^{2k}\delta_{1})\frac{a_{1}}{\delta_{1}}+(b+\lceil d\rceil
    2^{2k}\delta_{1})\left(\alpha-\frac{a_{1}}{\delta_{1}}\right)\\
    &=&(b+\lceil d\rceil
    2^{2k}\delta_{1})\frac{a_{1}}{\delta_{1}}+(b+d
    2^{2k}\delta_{1})\left(\alpha-\frac{a_{1}}{\delta_{1}}\right)+(\lceil
    d\rceil-d)(\delta_{1}\alpha-a_{1})\\
    &=&\mu-\frac{s}{\delta_{1}}+\lceil d\rceil
    2^{2k}a_{1}+\frac{s}{\delta_{1}}+(\lceil
    d\rceil-d)2^{2k}s|\delta_{1}\alpha-a_{1}|\\
    &=&\mu+\lceil d\rceil 2^{2k}a_{1}+(\lceil
    d\rceil-d)2^{2k}s|\delta_{1}\alpha-a_{1}|
  \end{eqnarray*}
  Using~(\ref{alphak}), we get
  \begin{eqnarray*}
    |\delta_{2}\alpha-(\mu+\lceil d\rceil 2^{2k}a_{1})|\in[0,2^{k}N^{-\frac{1}{2}}).
  \end{eqnarray*}
  Since $\lceil
  d\rceil<d+1=|\alpha-\tfrac{a_{1}}{\delta_{1}}|^{-1}2^{-2k}\delta_{1}^{-2}-b2^{-2k}\delta_{1}^{-1}+1$,
  (\ref{alphak})~yields the estimation $$\delta_{2}=b+\lceil d\rceil
  2^{2k}\delta_{1}<|\alpha-\tfrac{a_{1}}{\delta_{1}}|^{-1}\delta_{1}^{-1}+2^{2k}\delta_{1}\leq
  2^{k+1}N^{\frac{1}{2}}+2^{2k}\delta_{1}.$$ On the other hand
  $\delta_{2}\geq
  b+d2^{2k}\delta_{1}=|\alpha-\tfrac{a_{1}}{\delta_{1}}|^{-1}\delta_{1}^{-1}>2^{k}N^{\frac{1}{2}}$.
  Thus, $\delta_{2}\in M(b,k)\subseteq M_{\delta_{1}}(k)$. Set
  $L_{2}:=\lceil 2^{-k}\tfrac{N^{\frac{1}{2}}}{12}\rceil$. Then the
  set $E:=E_{\delta_{1},L_{1},\delta_{2},L_{2}}$ is an element of
  ${\mathcal E}_{3}$ and thus $E\in{\mathcal E}_{0}$.

  Before we can apply Lemma~\ref{cardinality}, we have to verify its
  conditions for the quadruple $(\delta_{1},L_{1},\delta_{2},L_{2})$.
  Since $(b,\delta_{1})=1$, also $(\delta_{1},\delta_{2})=1$.
  Moreover,
  \begin{eqnarray*}
    \delta_{2}\geq b+d\delta_{1}2^{2k}&=&|\alpha-\tfrac{a_{1}}{\delta_{1}}|^{-1}\delta_{1}^{-1}\\
    &>&(2^{-k}N^{-\frac{1}{2}}\delta_{1}^{-1})^{-1}\delta_{1}^{-1}\\
    &=&2^{k}N^{\frac{1}{2}}\\
    &>&\left\lceil\tfrac{2^{k}N^{\frac{1}{2}}}{12}\right\rceil=L_{1}.
  \end{eqnarray*}
  Thus, the conditions of Lemma~\ref{cardinality} are satisfied and
  the cardinality of the set $E:=\{j_{1}\delta_{1}+j_{2}\delta_{2}\mid
  j_{1}\in[0,L_{1}-1], j_{2}\in[0,L_{2}-1]\}$ can be estimated as
  follows: $|E|=L_{1}L_{2}\geq
  \tfrac{2^{k}N^{\frac{1}{2}}}{12}\tfrac{2^{-k}N^{\frac{1}{2}}}{12}=\tfrac{N}{144}$.
  Therefore, Lemma~\ref{fourier} proves
  $$|\widehat{\eins}_{-E}(\alpha)|\geq\frac{|E|}{2}\geq\frac{N}{288}.$$
\end{proof}


\begin{thebibliography}{99}

\bibitem[AS00]{ase} N.~Alon, J.~Spencer, and P.~Erd\H{o}s.  \newblock
  {\em The Probabilistic Method} (Second Edition).  \newblock John
  Wiley \& Sons, Inc., New York, 2000.

\bibitem[B81]{beck} J.~Beck.  \newblock Roth's Estimate of the
  Discrepancy of Integer Sequences is Nearly Sharp. {\em
    Combinatorica} {\bf 1}(4) (1981), 319-325.

\bibitem[BS95]{bs95} J. Beck and V. T. S\'os.  \newblock {\em
    Discrepancy theory}.  \newblock In R. Graham, M. Gr\"oschel and L.
  Lov\'asz, editors, {Handbook of Combinatorics}, 1405-1446. 1995.

\bibitem[DSW]{dsw} B.~Doerr, A.~Srivastav, P.~Wehr.  \newblock
  Discrepancies of Cartesian Products of Arithmetic Progressions. {\em
    Electronic Journal of Combinatorics}, to appear.

\bibitem[ES74]{es74} P.~Erd\H{o}s, J. Spencer.  \newblock {\em
    Probabilistic Methods in Combinatorics}.  \newblock Akademiai
  Kiado, Budapest, 1974.

\bibitem[H2004]{h04} N.~Hebbinghaus.  \newblock Discrepancy of Sums of
  Arithmetic Progressions. {\em Electronic Notes in Discrete
    Mathematics} {\bf 17C} (2004), pages 185-189.

\bibitem[Mat99]{matbook} J.~Matou\v{s}ek.  \newblock {\em Geometric
    Discrepancy}.  \newblock Springer, Berlin, 1999.

\bibitem[MS96]{ms96} J.~Matou\v{s}ek, J.~Spencer.  \newblock
  Discrepancy in Arithmetic Progressions. {\em Journal of the American
    Mathematical Society} {\bf 9}(1) (1996), 195-204.

\bibitem[P2006]{p06} A.~P\v{r}\'{i}v\v{e}tiv\'{y}.  \newblock
  Discrepancy of Sums of Three Arithmetic Progressions. {\em The
    Electronic Journal of Combinatorics} {\bf 13}(1) (2006), R5.

\bibitem[R64]{roth} K.~F.~Roth.  \newblock Remark Concerning Integer
  Sequences. {\em Acta Arithmetica} {\bf 9} (1964), 257-260.

\bibitem[V2002]{v02} B.~Valk\'o.  \newblock Discrepancy of Arithmetic
  Progressions in Higher Dimensions. {\em Journal of Number Theory}
  {\bf 92} (2002), 117-130.

\end{thebibliography}
\end{document}